\documentclass[oneside,11pt,reqno]{amsart}
\usepackage{amssymb, amscd}
\usepackage[all]{xy}
\usepackage{graphics}
\usepackage{textcomp}
 \theoremstyle{definition}
 
 \theoremstyle{remark}
 
\def\id{\mathrm{id}}

\def\lll{\mathbf{l}}
\def\nnn{\mathbf{n}}
\def\mmm{\mathbf{m}}
\def\coloneqq{:=}
\def\eqqcolon{=:}

\def\BZ{{\mathbb{Z}}}

\def\nn{\mbox{\boldmath ${n}$}}

\def\gg{\mbox{\boldmath ${g}$}}
\def\ff{\mbox{\boldmath ${f}$}}

\newlength{\vscaling} \newlength{\hscaling}
\usepackage{latexsym}
\usepackage{amssymb}

\def\FF{{\mathcal F}}
\def\GG{{\mathcal G}}

\def\P2{{P^{[2]}}}

\def\ll{{\mbox{\boldmath $\ell$}}}
\def\11{{\mbox{\boldmath $1$}}}

\def\eq{\begin{equation}}
\def\en{\end{equation}}

\def\eq#1\en{\begin{equation}#1\end{equation}}
\def\eqa#1\ena{\begin{eqnarray}#1\end{eqnarray}}

\def\qq{{\mbox{\boldmath $q$}}}

\newcommand{\mm}{{\mbox{\boldmath $m$}}}

\newcommand{\BB}{{\mathcal B}}

\newcommand{\QQ}{{\mathcal Q}}

\newcommand{\CC}{\mathcal{C}}

\newcommand{\RR}{\mathcal{R}}
\newcommand{\OO}{\mathcal{O}}
\newcommand{\PP}{\mathcal{P}}

\begin{document}
\begin{titlepage}
\vspace{1.cm}
\begin{flushright}
\end{flushright}
\vspace{1 cm}
\begin{centering}
\vspace{.43in} {\Large {Nonabelian bundle 2-gerbes}}\\
\vskip0.5cm
Branislav Jur\v co\\\vskip0.1cm
Max Planck Institute for Mathematics\\
Vivatsgasse 7, 53111 Bonn, Germany\\ \vspace{1.5in}
{\bf Abstract}\\
\end{centering}
\vspace{.1in}
 \noindent
We define 2-crossed module bundle 2-gerbes related to general Lie 2-crossed modules and discuss their properties. If $(L\to M\to N)$ is a Lie 2-crossed module and $Y\to X$ is a surjective submersion then an $(L\to M\to N)$-bundle 2-gerbe over $X$ is defined in terms of a so called $(L\to M\to N)$-bundle gerbe over the fibre product $Y^{[2]}=Y\times_X Y$, which is an $(L\to M)$-bundle gerbe over $Y^{[2]}$ equipped with a trivialization under the change of its structure crossed module from $L\to M$ to $1\to N$, and which is subject to further conditions on higher fibre products $Y^{[3]}$, $Y^{[4]}$ and $Y^{[5]}$. String structures can be described and classified using 2-crossed module bundle 2-gerbes.
\vspace{3.cm}

\end{titlepage}

\setcounter{page}{1}
\section*{Introduction}
The modest purpose of this paper is to introduce nonabelian bundle 2-gerbes related to 2-crossed modules \cite{Cond},  simultaneously generalizing abelian bundle 2-gerbes \cite{Stev2-gerbes}, \cite{StevPhD}, \cite{CareyMuWHighrBG} and crossed-module bundle gerbes \cite{AschieCantJur}, \cite{Jurco}. The idea is to describe objects in differential geometry, which would, in the terminology of \cite{BreenNotes},  correspond to the \v Cech cohomology classes in $H^1(X, L\to M\to N)$, i.e., the first \v Cech cohomology classes on a manifold $X$ with values in a Lie 2-crossed module $L\to M\to N$. What we want is a theory, which in the case of the 2-crossed module $U(1)\to 1\to 1$ reproduces the theory of abelian bundle 2-gerbes and in the case of a 2-crossed module $1\to M\to N$ reproduces the theory of crossed module bundle gerbes related to the crossed module $M\to N$ ($(M\to N)$-bundle gerbes). The latter requirement can slightly be generalized as follows. Let us assume a crossed module $L\stackrel{\partial}{\to} M$. If we put $A:=\ker \partial$ and $Q:=\rm coker \partial$ then we have a 4-term exact sequence of Lie groups $0\to A\to L\stackrel{\partial}{\to} M\to Q\to 1$ with abelian $A$. Let us assume that $A=U(1)$ is in the centre of $L$ and that the restriction to $U(1)$ of the action of $M$ on $L$ is trivial. Then we want that an $(U(1)\to L\to M)$-bundle 2-gerbe is the same thing as an $(L\to M)$ bundle gerbe twisted with an abelian bundle 2-gerbe \cite{AschJur}.

The paper is organized as follows.
In section 2, we briefly recall the relevant notions of a Lie crossed module and Lie 2-crossed module.
In section 3, relevant results on crossed module bundles and on crossed module bundle gerbes are collected. Let us mention that crossed module bundles are special kinds of bitorsors \cite{Gro}, \cite{Gi}, \cite{BreenBitors}, \cite{BreenNotes} and that crossed module bundle gerbes can be seen as a special case of gerbes with constant bands (this follows, e.g., from discussion in section 4.2. of \cite{BreenNotes} commenting on abelian bundle gerbes of \cite{MurrayBundleGerbes}, cocycle bitorsors of \cite{Ulbrich1}, \cite{Ulbrich2}, and bouquets of \cite{DuskinNonABinTopos}).
In section 4, 2-crossed module bundle gerbes are introduced as crossed module bundle bundles with an additional structure. 2-crossed module bundle gerbes are to 2-crossed module bundle 2-gerbes the same as crossed module bundles to crossed module bundle gerbes.
Finally, in section 5, 2-crossed module bundle 2-gerbes are introduced and their properties discussed, including their local description in terms of 3-cocycles similar to those of \cite{Dedecker3-dim} and \cite{BreenAst}, \cite{BreenNotes}. The example of a lifting bundle 2-gerbe is described in some detail. Also, we  discuss the relevance of 2-crossed module bundle 2-gerbes to string structures and their classification (see proposition \ref{Abelian4} and remark \ref{StringStr}). For the relevance of gerbes and abelian 2-gerbes to the string group and string structures see, e.g., \cite{BaezStevClass}, \cite{Bunke2}, \cite{BrylMcLau}, \cite{Jurco}, \cite{MurrayStevHigg}, \cite{StevStringGer} and \cite{Waldorf}. For discussions of abelian 2-gerbes in relation to quantum field theory and sting theory see, e.g., \cite{MickAnom}, \cite{CareyJohnsMurrStevWang}, \cite{CareyMuWHighrBG}, \cite{AschJur}.

Let us mention that in \cite{BreenAst} and \cite{BreenNotes} much more general 2-gerbes were introduced in the language of 2-stacks. These are generalizations of gerbes (defined as locally nonempty and locally connected stacks in groupoids \cite{Gi}, \cite{MoerdijkIntro}, \cite{BreenNotes}, \cite{BreenBitors}) and seem to be related rather to crossed squares than to 2-crossed modules. We hope to return to a discussion concerning a possible relation of our bundle 2-gerbes and the 2-gerbes of \cite{BreenAst} and \cite{BreenNotes} in the future. Also, we hope to discuss the relevant notion of a 2-bouquet elsewhere. Our task here is to describe nonabelian bundle 2-gerbes using a language very close to that of the classical reference books \cite{Koba-Nomizu}, \cite{Husemoller}. This will allow us introduce connection, curvature, curving etc. in the forthcoming paper \cite{JurcoGoogle2} using the language of differential geometry. For some further related work see, e.g.,  \cite{SatiSchrStash}, \cite{SSS}, \cite{SSS1}, \cite{RobSchr}, \cite{MarPick1}.

In this paper we work in the category of differentiable manifolds. In particular, all groups (with exception of the string group) are assumed to be Lie groups and all maps are assumed to be smooth maps. It would be to possible work with (for instance, paracompact Hausdorff) topological spaces, topological groups and continuous maps too. For this we would have to use a proper replacement of the notion of the surjective submersion $\wp:Y\to X$ in the definitions of crossed module bundle gerbes, 2-crossed module bundle gerbes and 2-crossed module bundle 2-gerbes. For instance, instead of surjective submersions we could consider surjective maps $\wp:Y\to X$ with a property that to each point $y\in Y$ there is a neighborhood $O$ of $\wp(y)$ with a section $\sigma: O\to Y$, such that $s(\wp(y))=y$. Such map may be called a surjective topological submersion. However, we should notice that there is also  another sort of map called topological submersion, which incidentally, looks like this: given $f:X \to Y$, for all $x \in X$ there is a neighborhood $U \subset X$ of $x$ such that when restricted to $U$ the restriction $f|U : U \to f(U)$   looks like a projection $U \times V \to U$ where $V$ is some (topological) vector space.
\vskip0.3cm
The present paper is based on my notes \cite{JurcoGoogle1}. It is a pleasure to thank MPIM for the opportunity to turn these into the present paper.
\vskip0.3cm
Further, it is a pleasure to thank Igor Bakovi\'c, David Roberts, Urs Schreiber, Danny Stevenson and Konrad Waldorf for discussions and comments.

\section{Crossed modules, 2-crossed modules}
Let us recall the notion of a crossed module of Lie groups (see, e.g., \cite{Brown},\cite{BrownR},\cite{PorterMenag}).
\subsection{Definition}\label{CM}
Let $L$ and $M$ be two  Lie groups. We say that $L$ is a crossed $M$-module if there
is a Lie group morphism $\partial_1: L\to M$ and a smooth action of $M$ on $L$
$(m,l)\mapsto \,^m\hskip-0cm l$ such that
$$^{\partial_1(l)}l' = ll'l^{-1}\,\, \mbox{(Peiffer condition)}$$ for $l,l'\in L$, and
$$\partial_1(^ml)= m \partial_1(l)m^{-1}$$ for $l\in L,
m\in M$ hold true. We will use the notation  $L\stackrel{\partial_1}{\to} M$ or $L\to M$ for the crossed module.

Let us also recall that a crossed
module is a special case of a pre-crossed module, in which the Peiffer condition doesn't necessarily hold.
There is an obvious notion of a morphism of crossed modules.
\subsection{Definition}\label{CMM} A morphism between crossed modules $L\stackrel{\partial_1}{\to} M$ and $L'\stackrel{\partial_1'}{\to} M'$
is a pair of Lie group morphisms $\lambda : L \to L'$ and $ \kappa : M \to M'$ such that the diagram
$$
\begin{CD}
 L@>
 {\partial_1}>> M\\
 @V {\lambda}
 VV @V VV{\scriptstyle\kappa}\\
L'@>
 {\partial_1'}>>M'\,
\end{CD}
$$
commutes, and for any $l \in L$ and $m \in M$ we have the following identity
\[ \lambda(^ml) = ^{\kappa(m)}\lambda(l). \]
\subsection{Remark}\label{2GROUP} A crossed module of Lie groups defines naturally a strict Lie 2-group. The set of objects is $C_0=\{*\}$, the set of 1-arrows is $C_1=M$ and the set of 2-arrows is $C_2= M\times L$. The "vertical" multiplication is given on $C_2$ by
$$(m,l_1)(\partial_1(l_1)m,l_2)=(m,l_1l_2)$$
and the "horizontal" multiplication is given by $$(m_1,l_1)(m_2, l_2) =(m_1m_2,l_1\,^{m_1}l_2)$$ See, e.g., \cite{BrowHigg} for more details on
the relation between crossed modules and strict Lie 2-groups.
\subsection{Definition} \label{2CM} The definition of a 2-crossed module of groups is according to Conduch\'e \cite{Cond}; see also, e.g., \cite{CondGeog}, \cite{MutPor}, \cite{BG}, \cite{PorterMenag}, \cite{RobSchr}).
A Lie 2-crossed module is a complex of Lie groups
\begin{equation}\label{2crossmod}
L\stackrel{\partial_1}{\to} M\stackrel{\partial_2}{\to} N
\end{equation}
together with smooth left actions by automorphism of $N$ on $L$ and $M$ (and on $N$ by conjugation), and the Peiffer lifting, which is an equivariant map $\{\,,\,\}: M\times M
\to L$ , i.e., $^n\{m_1, m_2\}=\{^n m_1, \,^n m_2\}$ such that:

i) (\ref{2crossmod}) is a complex of $N$-modules, i.e., $\partial_1$ and $\partial_2$ are $N$-equivariant and $\partial_1\partial_2 (l) =1$ for
$l\in L$,

ii) $m_1m_2m_1^{-1}=\partial_1\{m_1,m_2\}\,^{\partial_2(m_1)}m_2\eqqcolon\langle m_1, m_2\rangle$, for $m_1,m_2\in M$,

iii) $[l_1,l_2]\coloneqq l_1l_2l_1^{-1}l_2^{-1}=\{\partial_1 l_1, \partial_1 l_2\}$, for $l_1,l_2\in L$,

iv) $\{m_1 m_2, m_3\} = \{m_1, m_2 m_3 m_2^{-1}\} \,^{\partial_2(m_1)}\{m_2, m_3\}$, for $m_1, m_2, m_3 \in M$,

v) $\{m_1, m_2 m_3\} = \{m_1, m_3\}\,^{m_1 m_2 m_1^{-1}}\{m_1, m_3\}$, for $m_1, m_2, m_3 \in M$,

vi) $\{\partial_1(l), m\}\{m,
\partial_1(l)\}= l^{\partial_2(m)}(l^{-1})$, for $m\in M, l\in L$, and

\noindent wherein the notation $^nm$ and $^nl$ for left actions of the element $n\in N$ on elements $m\in M$ and $l\in L$ has been used.
Also, also let us notice that $^ml \coloneqq l\{\partial_1(l)^{-1}, m\}$ defines a left action of $M$ on $L$ by automorphisms. This is a consequence of the
other axioms and is proved in \cite{Cond}, \cite{BG}, where it is also shown that, equipped with this action, $L\stackrel{\partial_1}{\to} M$ defines
a crossed module.
\subsection{Remark} \label{2ndCM}In addition to the crossed module $L\stackrel{\partial_1}{\to} M$, there is an another crossed module that can be associated with the 2-crossed module $L\stackrel{\partial_1}{\to} M\stackrel{\partial_2}{\to} N$. By definition, we see that  $M\stackrel{\partial_2}{\to} N$ is a (special) pre-crossed module in which the Peiffer condition is satisfied only up to the Peiffer lifting. Hence,
$\partial_1(L)\backslash M \stackrel{\partial_2'} {\to} N$, with the induced Lie group homomorphism $\partial_2'$ and with the induced action of $N$ on $\partial_1(L)\backslash M$, is a crossed module.

There is an obvious notion of a morphism of 2-crossed modules.
\subsection{Definition}\label{2CMM} A morphism between 2-crossed modules $L\stackrel{\partial_1}{\to} M\stackrel{\partial_2}{\to} N$ and $L'\stackrel{\partial_1'}{\to} M'\stackrel{\partial_2'}{\to}N'$ is a triple of Lie group morphisms $L\to L'$, $M\to M'$ and $N\to N'$making up, together with the maps $\partial_1$, $\partial_1'$, $\partial_2$ and $\partial_2'$ a commutative diagram
\begin{equation}\label{mor2crossmod}
\begin{CD}
 L@>
 {\partial_1}>> M@>
 {\partial_2}>> N\\
 @V {\lambda}
 VV @V \mu VV @V VV{\scriptstyle\nu} \\
L'@>
 {\partial_1'}>>M'@>
 {\partial_2'}>>N'\,
\end{CD}
\end{equation}
and being compatible with the actions of $N$ on $M$ and $L$ and of $N'$ on $M'$ and $L'$, respectively and with the respective Peiffer liftings.
\subsection{Remark}\label{Gray} A 2-crossed module of Lie groups defines naturally a Gray 3-groupoid with a single object. Hence, the set of objects of the Gray 3-groupoid is $C_0=\{*\}$. The set of 1-arrows  is $C_1=N$, the set of 2-arrows is $C_2= N\times M$ and the set of 3-arrows is $C_3= N\times M\times L$. For the construction and for more details on the relation between 2-crossed modules and Gray 3-groupoids see \cite{KampsPorter}, \cite{BG}, \cite{MutPor}, \cite{MarPick1}.
There are two "vertical" multiplications and one "horizontal" multiplication. The vertical multiplications are determined by the crossed module $L\to M$. On $C_3$, the two vertical multiplications are given by
$$(n,m,l_1)(n,\partial_1(l_1)m,l_2)=(n,m,l_1l_2)$$
and
$$(n,m_1,l_1)(\partial_2(m_1)n, m_2, l_2) =(n,
m_1m_2,l_1\,^{m_1}l_2)$$ and the horizontal multiplication is given by
$$(n_1,m_1,l_1)(n_2,m_2,l_2)=(n_1n_2, m_1\,^{n_1}m_2,
l_1\,^{m_1}(\,^{n_1}l_2))$$
\section{Crossed module bundle gerbes}
Let $X$ be a (smooth) manifold. Crossed module bundle gerbes have been introduced, for instance, in \cite{Jurco}, \cite{AschieCantJur}. These can be seen as generalizations of abelian bundle gerbes \cite{MurrayBundleGerbes}, \cite{MurrayIntr}.
If $(L\stackrel{\partial_1}{\to} M)$ is a crossed module of  Lie groups, $X$ a manifold and $P\to X$ a left principal $L$-bundle, we can change the structure group of $P$ from $L$ to $M$, in order to obtain a left principal $M$-bundle $P'= M\times_{\partial_1}P$ defined as follows. Points $p' \in P'$ correspond  to  equivalence classes $[m, p] \in M\times_{\partial_1}P$ with the equivalence relation on $M\times P$ given by $(m, p)\sim (m\partial_1(l), l^{-1}p)$. Obviously, the principal left $M$-action is given by $M\times P' \to P'$, $m'\times [m,p]\mapsto [m'm, p]$.
\subsection{Definition}\label{CMB}
Let $(L\stackrel{\partial_1}{\to} M)$ be a crossed module of  Lie groups and $X$ a manifold. Let $P\to X$ be left principal $L$-bundle, such that the principal $M$-bundle $M\times_{\partial_1}P$ is trivial with a  trivialization defined by a section (i.e., a left $L$-equivariant smooth map) $\mm: P\to M$. We call the couple $(P,\mm)$ an
$(L\to M)$-bundle.
\subsection{Remark} \label{GroupoidBundle}If we think about the crossed module $L\to M$ as about a groupoid with the set of objects $M$ and the set of arrows $M\times L$ then a crossed module bundle is the same thing as a principal groupoid bundle.
\subsection{Definition}\label{Iso1}Two $(L\to M)$-bundles $(P,\mm)$ and $(P',\mm')$ over $X$  are isomorphic if they are isomorphic as left $L$-bundles by an isomorphism $\ll: P\to P'$ and $\mm'\ll = \mm$. An $(L\to M)$-bundle is trivial if it is isomorphic to the trivial $(L\to M)$-bundle $(X\times L, \partial_1{\rm pr}_L)$.
\subsection{Example}\label{NLT1}Let us notice, that a general $(L\to M)$-bundle is not necessarily locally trivial, although it is locally trivial as a left principal $L$-bundle. For instance, for a function $m: X\to M$ such that ${\rm Im}(m)$ is not a subset of ${\rm Im}(\partial_1)$ the $(L\to M)$-bundle $(X\times L, \partial_1{\rm pr}_L.m{\rm pr}_X)$ is locally non-trivial. We will refer to such an  $(L\to M)$-bundle as to an $(L\to M)$-bundle defined by the $M$-valued function $m$. Two such $(L\to M)$-bundles are isomorphic iff their respective functions $m$ and $m'$ are related by an $L$ valued function $l$ on $X$ by $m'=\partial_1(l)m$.
\subsection{Example}\label{1toG} A $(1\to G)$-bundle is the same thing as a $G$-valued function.
\subsection{Example}\label{LT1}A couple $(T,\lll)$, where $T$ is a trivial left principal $L$-bundle and $\lll: T\to L$ its trivializing section, defines an $(L\to M)$-bundle with the section $\mm=\partial_1\lll: T\to M$. $(T,\lll)$ is a trivial $(L\to M)$-bundle.
\subsection{Example}\label{NSCMB} Let $L$ be a normal subgroup of $M$. The adjoint action of $M$ restricted to $L$ defines a crossed module structure on $L\to M$ with $\ker \partial_1=1$. Let $L$ be also a closed subgroup of $M$. We put $G\coloneqq M\backslash L$, so that we have an exact sequence of Lie groups $1\to L\to M\stackrel{\pi}{\to} G\to 1$, such that  $M$ is a left principal $L$-bundle over $G$. Moreover, $(M\to G, \mm)$ with $\mm = {\rm{id}}_M$ is an $(L\to M)$-bundle.
\subsection{Pullback} \label{Pull1}Obviously, a pullback of an $(L\to M)$-bundle is again an  $(L\to M)$-bundle. Pullbacks preserve isomorphisms of crossed module bundles, in particular a pullback of a trivial $(L\to M)$-bundle is again a trivial $(L\to M)$-bundle.
\subsection{Change of the structure crossed module} \label{Change1}If $(L\to M) \to (L'\to M')$ is a morphism of crossed modules, there is an obvious way to construct, starting from an $(L\to M)$-bundle $(P,\mm)$, an $(L'\to M')$-bundle $(L'\times_\lambda P,\kappa\mm)$ where $\lambda: L\to L'$ and $\kappa: M\to M'$ define the morphism of the two crossed modules. Obviously, the change of the structure crossed module preserves isomorphisms of crossed module bundles.
\subsection{1-cocycles} \label{1-cocycle}Consider an $(L\to M)$-bundle $(P,\mm)$ and  a trivializing covering $\coprod {O}_{i}=X$ of the left principal $L$-bundle $P$. Let $\sigma_i: P|O_i \to L$ be the trivializing sections of $L$ and $l_{ij}=\sigma_i^{-1} \sigma_j: O_{i}\cap O_{j}\to L$ be the corresponding transition functions. We put $m_i=\partial_1(\sigma_i)^{-1}\mm$, which obviously gives an $L$-valued function on $O_i$. We have $\partial_1 (l_{ij})=m_i m_j^{-1}$. Hence the $(L\to M)$-bundle $(P,\mm)$ can be described by a 1-cocycle given by transition functions $(m_i,l_{ij})$, $m_i: O_i \to M$, $l_{ij}: O_{ij}=O_{i}\cap O_{j}\to L$ satisfying on nonempty $O_{ij}=O_{i}\cap O_{j}$
$$\partial_1(l_{ij})=m_i m_j^{-1}$$ and on nonempty $O_{ijk}=O_{i}\cap O_{j}\cap O_{k}$
$$l_{ij}l_{jk}=l_{ik}$$
Transition functions $ (m_i,l_{ij})$ and $(m'_i,l'_{ij})$ corresponding to two isomorphic $(L\to M)$-bundles are related by
$$m'_i = \partial_1(l_i)m_i$$
and
$$l'_{ij} = l_i l_{ij}l_i^{-1}$$
A trivial $(L\to M)$-bundle is described by transition functions $(\partial_1(l_i), l_i l_j^{-1})$.

On the other hand, given transition functions $(m_i,l_{ij})$ we can reconstruct an $(L\to M)$-bundle. We define a left principal $L$-bundle $P$ with the total space formed by equivalence classes of triples $[x,l,i]$ with $x\in O_i$, $l\in L$ under the equivalence relation $(x,l,i)\sim (x',l',j)$ iff $x=x'$ and $l'=ll_{ij}$. The principal left $L$-action is given by $l'[x,l,i]=[x,l'l,i]$. Now we put $\mm([x,l,i])=\partial_1 (l) m_i(x)$. $(P,\mm)$ is an $(L\to M)$-bundle.

With the two above constructions it is not difficult to prove that  the category of $(L\to M)$-bundles is equivalent to the category of 1-cocycles.

\subsection{Lifting crossed module bundle}\label{Lifting1} Let, as above in (\ref{NSCMB}), $L$ be a normal closed subgroup of $M$ and $L\to M$ the corresponding crossed module. Consider a $G$-valued function $g: X\to G$. The pullback $g^*(M,{\rm{id}}_M)$ of the $(L\to M)$-bundle $\pi:M\to G$ is an $(L\to M)$-bundle on $X$. It is the obstruction to a lifting of the $G$-valued function $g$ to some $M$-valued function. To go in the opposite direction, we notice that we have an obvious morphism of crossed modules $(L\to M)\to (1\to G)$. Under the change  of the structure crosse module of an $(L\to M)$-bundle $(P,\mm)$ from $(L\to M)$ to $(1\to G)$, the section $\mm$ becomes an $L$-invariant $G$-valued function $\pi\mm$ on $P$. Hence, it is identified with a $G$-valued function $g$ on $X$. Two isomorphic $(L\to M)$-bundles give the same $G$-valued function. The two constructions are inverse to each other up to an isomorphism of $(L\to M)$-bundles.

It is now easy to give a local description of lifting crossed module bundles.
Let $\{O_i\}_i$ be a covering on $X$. Let $P$ be a $(L\to M)$-bundle described by transition functions $(l_{ij}, m_i)$. Then, since $\pi\partial_1 =1$, we have $\pi(m_i)=\pi(m_j)$ Hence, the collection of local functions $\{\pi(m_i)\}_i$ defines a $G$-valued function $X$. To go in the opposite direction, let $g$ be $G$-valued function on $X$ . Let $O_i$ be the trivializing covering of the pullback  principal bundle $L$-bundle $g^*(M)$. The function $g$ can now be described by a collection of local functions $g_i: O_i\to G$ such that $g_i=g_j$ on $O_{ij}$. Hence, we have local functions $m_i: O_i\to M$, the local sections of $g^*(M)$ such that $\pi(m_i)=g_i$, which are related on double intersections $O_{ij}$ as $m_i=\partial_1(l_{ij})m_j$ by $L$-valued functions $l_{ij}:O_{ij}\to L$, the transition functions of the principal bundle $M\stackrel{\pi}{\to}G$ fulfilling the 1-cocycle condition $l_{ij}l_{jk}=l_{ik}$ on $O_{ijk}$.

Concerning crossed module bundles, we have the following lemma and proposition (see \cite{AschieCantJur}).
\subsection{Lemma} \label{RightAction1}The $(L\to M)$-bundle $(P,\mm)$ is also a right principal $L$-bundle with the right action of $L$ given by $p.l = \, ^{{\scriptstyle{\mm}}(p)}(l).p$ for $p\in P, l\in L$. The left and right actions commute. The section $\mm$ is $L$-biequivariant.
\subsection{Proposition}\label{mult1} Let $(P, \mm)$ and $(\tilde P,\tilde \mm)$ are two $(L\to M)$-bundles over $X$. Let us define an equivalence relation on the Whitney sum  $P\oplus \tilde P=P \times_X \tilde P$ by $(pl, \tilde p)\sim (p, l\tilde p)$, for $(p,\tilde p)\in P\oplus \tilde P$ and $l\in L$. Then $(P\tilde P\coloneqq(P\oplus \tilde P)/\sim, \mm\tilde \mm)$ with $\mm\tilde \mm([p,\tilde p])\coloneqq\mm(p)\tilde \mm(\tilde p)$ is an $(L\to M)$-bundle.

\subsection{Remark}\label{Group1} The set of isomorphism classes of $(L\to M)$-bundles equipped with the above defined product is a group. The unit is given by the class of the trivial bundle $(X\times L, \partial_1\rm{pr}_L)$. The inverse is given by the class of $(L\to M)$-bundle $(P^{-1},\mm^{-1})$ with $P^{-1}$ having the same  total space as $P$, the left $L$-action on $P^{-1}$ being the inverse of the right $L$-action on $P$ and the trivializing section $\mm^{-1}$ being the composition of the inverse in $M$ with the trivializing section $\mm$. Let us notice that in the case of an exact sequence $1\to L  \to M \to N\to 1$ as above (\ref{NSCMB}) this group structure is compatible with the group structure of $G=L\backslash M$-valued functions  with point-wise multiplication.
\subsection{Example} \label{ProdNLT1}If $\PP=(P=X\times L, \partial_1{\rm pr}_L.m{\rm pr}_X)$ and $\PP'=(P'=X\times L, \partial_1{\rm pr}_L.m'{\rm pr}_X)$ are $(L\to M)$-bundles defined by two respective $M$-valued functions $m$ and $m'$ on $X$ (\ref{NLT1}) then the product $\PP\PP'$ is explicitly described again as an $(L\to M)$-bundle defined by the function $mm'$ by identifying $[(x,l), (x,l')]\in PP'$ with $(x, l\,^ml')\in X\times L$.
\subsection{Product on 1-cocycles}\label{Prod1-Cocyles}
Transition functions $(\bar m_i,\bar l_{ij})$ of the product of two $(L\to M)$-bundles described by transition functions $(m_i,l_{ij})$  and $(\tilde m_i,\tilde l_{ij})$ are given by
$${\bar m_i}={m_i}{\tilde m_i}$$
and
$$\bar l_{ij}=l_{ij}\,^{m_i}\tilde l_{ij}$$
Transition functions of the inverse crossed module bundle are $(^{m_j^{-1}}l_{ij}^{-1}=\,^{m_i^{-1}}l_{ij}^{-1}, m_i^{-1}).$

\subsection{1-cocycles as functors}\label{Funct1}The crossed module $(L\to M)$ defines naturally a topological category (groupoid) $\CC$ with the set of objects $C_0=L$ and the set of arrows $C_1=M\times L$. Let us consider the topological category $\OO$ (groupoid) defined by the good covering $\{O_i\}$ of $X$ with objects $x_i\coloneqq(x, i|\,x\in O_i)$ and exactly one arrow from $x_i$ to $y_j$ iff $x=y$. Then a 1-cocycle is the same thing as a  continuous functor from $\OO$ to $\CC$.
Further, if $2\BB$ is a strict topological 2-category, then the category of 2-arrows with the vertical composition is naturally a topological category $\BB$. The horizontal composition in  $2\BB$ defines a continuous functor from the cartesian product $\BB \times \BB$ to $\BB$. Thus, in case $\BB=\CC$ it defines naturally a multiplication on functors $\OO \to \CC$ (i.e., on transition functions), which is the same as the one defined above (\ref{mult1}).
\subsection{Crossed module bundle gerbes}\label{CMBG} Let $Y$ be a manifold. Consider a  surjective submersion $\wp~:~Y\rightarrow X$, which in particular admits local sections. Let $\{{O}_{i}\}$ be the corresponding covering of $X$
with local sections $\sigma_{i}:{O}_i\rightarrow Y$,
i.e., $\wp\sigma_i=id$.
We also consider $Y^{[n]}=Y\times_X Y\times_X Y\ldots\times_X Y$, the n-fold
fibre product of $Y$, i.e.,
$Y^{[n]}\coloneqq\{(y_1,\ldots y_n)\in Y^n\;|\;\wp(y_1)=\wp(y_2)=
\ldots\wp(y_n)\}$.
Given an  $(L\to M)$-bundle $\PP=(P,\mm)$ over $Y^{[2]}$
we denote by ${\PP}_{12}=p_{12}^*(\PP)$
the crossed module bundle on $Y^{[3]}$ obtained as a pullback
of  $\PP$ under $p_{12}:Y^{[3]}\rightarrow Y^{[2]}$ ($p_{12}$ is the identity on
its first two arguments);
similarly for ${\PP}_{13}$ and ${\PP}_{23}$.
Consider a quadruple
$({\PP},Y,X,\ll)$, where ${\PP}=(P,\mm)$ is a crossed module bundle, $Y\to X$ a surjective submersion
and $\ll$ an isomorphism of
crossed module bundles
$\ll :  {\PP}_{12}{\PP}_{23}\to{\PP}_{13}$.
We now consider bundles ${\PP}_{12},\,{\PP}_{23},\,{\PP}_{13},\,{\PP}_{24},
\,{\PP}_{34},\,{\PP}_{14}$
on $Y^{[4]}$
relative to the projections $p_{12}~:~Y^{[4]}\rightarrow Y^{[2]}$ etc.
and also the crossed module isomorphisms
$\ll_{123}, \,\ll_{124},\, \ll_{123}, \,\ll_{234}$ induced by projections
$p_{123}:Y^{[4]}\rightarrow Y^{[3]}$ etc.
\subsection{Definition}\label{DefCMBG}
The quadruple
$({\PP},Y,X,\ll)$, where $Y\to X$ is a surjective submersion, ${\PP}$ is a crossed module bundle over $Y^{[2]}$, and $\ll :  {\PP}_{12}{\PP}_{23}\to{\PP}_{13}$ an isomorphism of
crossed module bundles over $Y^{[3]}$, is called a crossed module bundle gerbe if $\ll$ satisfies the
cocycle condition (associativity) on $Y^{[4]}$
\begin{center}
\begin{equation}\label{crossmodbunger}
\begin{CD}
 \PP_{12}\PP_{23}\PP_{34}@>
 {\ll_{234}}>> \PP_{12}\PP_{24}\\
 @V {\ll_{123}}
 VV @V  VV \hskip-0.3cm{\scriptstyle{\ll_{124}}}\\
\PP_{13}\PP_{34}@>
 {\ll_{134}}>>\PP_{14}
\end{CD}
\end{equation}
\end{center}
\subsection{Abelian bundle gerbes}\label{Ab}Abelian bundle gerbes as introduced in \cite{MurrayBundleGerbes}, \cite{MurrayIntr} are $(U(1)\to 1)$-bundle gerbes. More generally, if $A\to 1$ is a crossed module then $A$ is necessarily
an abelian group and an abelian bundle gerbe can be identified as an $(A\to 1)$-bundle gerbe.
\subsection{Example} \label{1toG2} A $(1\to G)$-bundle gerbe is the same thing as a $G$-valued function $g$ on $Y^{[2]}$ (\ref{1toG}) satisfying on $Y^{[3]}$ the cocycle relation $g_{12}g_{23}= g_{23}$ and hence, a principal $G$-bundle on $X$ (more precisely a descent datum of a principal $G$-bundle).
\subsection{Pullback} If $f :X\to X'$ is a map then we can pullback $Y\to X$ to $f^*(Y)\to X'$ with a map $\tilde f: f^*(Y)\to Y$ covering $f$.
There are induced maps $\tilde f^{[n]} : f^*(Y )^{[n]} \to Y^{[n]}$.
Then the pullback $f^*({\PP},Y,X,\ll)\coloneqq (\tilde f^{[2]*}{\PP},f^*(Y),f(X),\tilde f^{[3]*}\ll)$ is again an $(L\to M)$-bundle gerbe.
\subsection{Definition}\label{Iso2}
Two crossed module bundle gerbes $(\PP,Y,X,\ll)$ and $(\PP',Y',X,\ll')$ are
stably isomorphic if there exists a crossed module bundle $\QQ \to \bar Y=Y\times_X Y'$ such that over $\bar Y^{[2]}$ the crossed module bundles $q^*\PP $ and
$\QQ_1 q'^*\PP'\QQ_2^{-1}$ are isomorphic. The corresponding isomorphism
$\tilde\ll: q^*\PP \to \QQ_1 q'^*\PP'\QQ_2^{-1}$ should satisfy on $\bar Y^{[3]}$ (with an obvious abuse of notation) the condition
$$\tilde\ll_{13}  \ll = \ll' \tilde\ll_{23} \tilde\ll_{12}$$
Here $q$ and $q'$ are projections onto first and second factor of $\bar Y=Y\times_X
Y'$ and $\QQ_1$ and $\QQ_2$ are the pullbacks of  $\QQ \to \bar Y$ to $\bar Y^{[2]}$ under respective
projections form $\bar Y^{[2]}$ to $\bar Y$ etc.

A crossed module bundle gerbe $(\PP,Y,X,\ll)$ is called trivial if it is stably isomorphic to the trivial crossed module bundle gerbe $((Y^{[2]}\times L,\partial_1{\rm pr}_L), Y, X, 1)$. Pullbacks preserve stable isomorphisms, in particular a pullback of a trivial crossed module bundle gerbe is again a trivial crossed module bundle gerbe.
If $Y=X$ then the crossed module bundle gerbe is trivial.
\subsection{Definition}\label{IsoIso} Let $({\PP},Y,X,\ll)$ and $({\PP'},Y',X,\ll')$ be two crossed module bundle gerbes and $(\QQ, \tilde\ll_\QQ)$ and $(\RR, \tilde\ll_\RR)$ two stable isomorphisms between them. We call $(\QQ, \tilde\ll_\QQ)$ and $(\RR, \tilde\ll_\RR)$ isomorphic if there is an isomorphism $\underline\ll:\QQ \to \RR$ of crossed module bundles on $\bar Y=Y\times_X Y'$ such that (with an obvious abuse of notation) the diagram
\begin{equation}\label{isocrossmodbunger}
\begin{CD}
 q^*\PP@>
 {\tilde\ll_\QQ}>> \QQ_1{q'}^*\PP'\QQ_2^{-1}\\
 @V{\rm id}
 VV @V  VV \hskip-0.5cm{\scriptstyle{\underline\ll_1}{\underline\ll_2}^{-1}}\\
q^*\PP@>
 {\tilde\ll_\RR}>>\RR_1{q'}^*\PP'\RR_2^{-1}
\end{CD}
\end{equation}
is commutative.
\subsection{Remark}\label{Indep1}  Let $\wp: Y'\to X$ be an another surjective submersion and $f: Y'\to Y$ a map such that $\wp'=\wp f$. Then crossed module bundle gerbes $\GG_f=(f^*\PP, Y', X, f^{[3]*}\ll)$ and $\GG=(\PP, Y, X, \ll)$ are stably isomorphic. This can be easily seen by noticing first that $\GG$ is stably isomorphic to itself and then using the obvious fact that pullbacks of crossed module bundles commute with their products \cite{AschieCantJur}. It follows that locally each crossed module bundle gerbe $\GG$ is trivial. For this, take a point $x\in X$ and its neighborhood $O\subset X$ such that there exists a local section $\sigma:O\to Y$ of $\wp$. Over $O$ we have the bundle gerbe $\GG_O$, the restriction of $\GG$ to $O$. Now we can put $Y'\coloneqq O$ and $\wp' := \id_O$ and we have $\wp\sigma=\wp'$. It follows that $\GG_\sigma$ is stably isomorphic to $\GG_O$. However $\GG_\sigma$ is trivial, because of $Y'=O$.
\subsection{Change of the structure crossed module} \label{Change2}If $(L\to M) \to (L'\to M')$ is a morphism of crossed modules, there is an obvious way to construct starting from an $(L\to M)$-bundle gerbe an $(L'\to M')$-bundle gerbe by changing the structure crossed module of the corresponding $(L\to M)$-bundle over $Y^{[2]}$. Obviously, the change of the structure crossed module preserves stable isomorphisms of crossed module bundle gerbes.

\subsection{2-cocycles}\label{2-cocycle}
Locally, bundle gerbes can be described in terms of 2-cocycles as follows.
First, let us notice that the trivializing cover $\{O_i \}$ of the map $\wp: Y\to X$ defines a new surjective submersion $\wp': Y'=\coprod O_i\to X$. The local sections of $Y\to X$ define a map $f:Y'\to Y$, which is compatible with the maps $\wp$ and $\wp'$, i.e., such that $\wp f = \wp'$. We know that crossed module bundle gerbes $\GG_f$ and $\GG$ are stably isomorphic. Hence, we can again assume $Y=\coprod O_i$. For simplicity, we assume that the covering $\{O_i \}$ is a good one. Then the crossed module bundle gerbe can be described by a 2-cocycle $(m_{i j}, l_{i j
k})$ where the  maps $m_{i j}: O_{ij}
\to M$ and $l_{i j k}: O_{ij
k} \to L$ fulfill the following conditions
$$m_{i j}m_{j k}= \partial_1 (l_{i j k})
m_{i k} \hskip0.4cm {\rm on}\hskip0.4cm O_{ijk}
$$
and
$$l_{i j k}l_{i k l} = \,^{m_{i
j}}\hskip-0.05cm
l_{j k l} l_{i j l} \hskip0.4cm {\rm on}\hskip0.4cm
O_{ijkl}$$
Two crossed module bundle gerbes are stably isomorphic if
their respective 2-cocycles
 $(m_{i j}, l_{i j k})$ and
 $(m_{i j}', l_{i j k}')$ are related
 by
 $$m_{i j}'=m_{i}\partial_1(l_{i j})
 m_{i j}m_{j}^{-1}$$
 and
 $$ l_{i j k}'= \,^{m_i}l_{i j}
 \,^{m_im_{ij}}l_{j k}
 \,^{m_i}l_{i j k}\,^{m_i}l_{i k}^{-1}$$
 with $m_i : O_i \to M$ and $l_{ij} : O_{ij}
 \to L$.
A trivial crossed module bundle gerbe is described by transition functions
$$m_{i j}=m_i \partial_1(l_{i j}) m_j^{-1}$$ and $$ l_{i j k}=\,^{m_i}l_{i j}
 \,^{m_i}l_{j k}
 \,^{m_i}l_{i k}^{-1}$$
Two collections of stable isomorphism data $(m_i, l_{ij})$ and $(m_i', l_{ij}')$ are isomorphic if
$$m_i' = \partial_1 (l_i)m_i $$
$$l_{ij}'= l_i l_{ij} l_{j}^{-1}$$
Now we briefly describe how an $(L\to M)$-bundle gerbe can be reconstructed from transition functions $(m_{ij}, l_{ijk})$. Put $Y=\coprod O_i$. On each nonempty $O_{ij}$ consider the $(L\to M)$-bundle $\PP_{ij}$  defined by the
function $m_{ij}:O_{ij}\to M$ as in (\ref{NLT1}). Hence, on $Y^{[2]}$ we have the $(L\to M)$-bundle given by $\PP=\coprod_{ij}\PP_{ij}$. Now we recall the explicit descriptions of the multiplication (\ref{ProdNLT1}) and isomorphisms (\ref{NLT1}) of two $(L\to M)$-bundles defined by their respective $M$-valued functions. Using the 2-cocycle relations, it is now straightforward to show that the collection of functions $l_{ijk}$ defines an isomorphism of $\PP_{12}\PP_{23}$ and $\PP_{13}$ on $Y^{[3]}$ satisfying the associativity condition on $Y^{[4]}$ (compare, e.g., to Theorem 3.1 in \cite{MoerdijkIntro})

Also, it is straightforward to check that stable isomorphism classes of $(L\to M)$-bundle gerbes are one to one with stable isomorphism classes of transition functions, the isomorphism being given by the two above described constructions.
Actually, when considering isomorphisms of stable isomorphisms, we have the respective 2-categories of $(L\to M)$-bundle gerbes and transition functions. The correspondence between $(L\to M)$-bundle gerbes and the transition function can be formulated in the framework of 2-categories similarly to \cite{BreenAst}, but we will not do this here.
Further, if we consider the topological $\OO$ category $\OO$ defined by the good covering $\{O_i\}$ of $X$.
Then a 2-valued cocycle can be seen as a continuous normal pseudo-functor from $\OO$ to the bicategory defined by the crossed module $L\to M$.
\subsection{Lifting crossed module bundle gerbe}\label{Lifting2}
Let $L\to M$ be a crossed module associated with a closed normal subgroup $L$ of $M$ (\ref{Lifting1}). We have a Lie group extension
$$1\to L\stackrel{\partial_1}{\to} M\stackrel{\pi}{\to} G \to 1$$
and also the $(L\to M)$-bundle $M\stackrel{\pi}{\to} G$.
Let $E\to X$ be a (locally trivial) left principal $G$-bundle over $X$.  As a principal $G$-bundle $E$ defines a (division) map $g:E^{[2]} \to G$ which gives for two elements in a fibre of $P$ the group element relating them. The pullback $\PP =g^*(M, {\rm{id}}_M)$ of the $(L\to M)$-bundle $M\to G$ gives an $(L\to M)$-bundle on $E^{[2]}$; $\PP$ is by definition the lifting $(L\to M)$-bundle corresponding to the division map $g$ (\ref{Lifting1}). It follows that the crossed module bundles $\PP_{12}\PP_{23}$ and $\PP_{13}$ are isomorphic on $Y^{[3]}$. This follows from the above mentioned fact that, in case of Lie groups $L$, $M$ and $G$ as above, isomorphism classes of $(L\to M)$-bundles are one to one to $G$-valued functions and that this correspondence respects the respective multiplications (\ref{Lifting1}). Such an isomorphism $\ll$  fulfills the associativity condition because of $\ker(\partial_1)=1$. Hence, we have a crossed module bundle gerbe $\GG$, which can be seen as an obstruction to a lifting of the principal $G$ bundle $E$ to some principal $M$-bundle. Also, it is easily seen that lifting two isomorphic $G$-bundle gerbes leads to stably isomorphic $(L\to M)$ bundle gerbes.  On the other hand, if we have a crossed module $L\to M$ with a trivial kernel of $\partial_1$ and hence fitting the exact sequence with $G={\rm coker}{\partial_1}$ we can change the structure crossed module from $L\to M$ to $1\to G$ in a crossed module bundle gerbe in order to get a principal $G$-bundle on $X$. These two constructions are inverse to each other on sets of stable isomorphism classes of $(L\to M)$-bundle gerbes (with $(L\to M)$  as above) and isomorphism classes of principal $G$-bundles. Finally, given three principal $G$-bundles $E, E''$ and $E'''$ and isomorphisms $E\stackrel{f}{\to} E'$, $E'\stackrel{f'}{\to} E''$ and $E''\stackrel{f''}{\to} E'''$ such that $f'f=f''$ the corresponding lifting crossed module bundle gerbes $\GG, \GG''$ and $\GG'''$ will be stably isomorphic, but the respective stable isomorphisms $\ff\ff'$ and $\ff''$ will be only isomorphic in general.

It is also easy to give a local description of lifting crossed module bundle gerbes. Let again $\{O_i\}_i$ be a good covering of $X$. Let us consider an $(L\to M)$-bundle gerbe described by transition functions $(m_{i j}, l_{i j k})$. Then $\pi(m_{i j})\pi(m_{jk})= \pi(m_{ik})$. Hence, we have a principal  $G$-bundle with transition functions $g_{i j}=\pi(m_{i j})$. To go in the opposite direction, let $g_{i j}$ be transition functions of a principal $G$-bundle.  Since the double intersections $O_{ij}$ are contractible we can choose lifts $m_{i j}$ of transition functions $g_{i j}$.  On $O_{ijk}$ these will be related by $m_{i j}m_{jk}=\partial_1 (l_{ijk})m_{i k}$ with $L$-valued functions $l_{ijk}$  which, because of $\ker \partial_1 =1$, necessarily satisfy  the requested compatibility condition on $O_{ijkl}$ (\ref {CMBG}).

\section{2-crossed module bundle gerbes}
Let $(L\to M\to N)$ be a Lie 2-crossed module and $\GG$ be an $(L\to M)$-bundle gerbe over
$X$. From the definition of the 2-crossed module we see immediately that maps $L\to 1$ and $\partial_2: M \to N$ define a morphisms of crossed modules $\mu: (L\stackrel{\partial_1}{\to} M)\to (1\to N)$. Thus, we have the following trivial lemma (\ref{1toG}):
\subsection{Lemma}\label{lemma1}
$\mu(\GG)$ is a principal $N$-bundle on $X$. If $\GG$ and $\GG'$ are stably isomorphic, then $\mu(\GG)$ and $\mu(\GG')$ are isomorphic.
\subsection{Definition} \label{2CMBG}Let $\GG$  be an $L\to M$-bundle gerbe such that the principal bundle $\mu(\GG)$ over $X$ is trivial with a section $\nn: \mu(\GG)\to N$. We call the pair $(\GG, \nn)$ a 2-crossed module bundle gerbe.
\subsection{Remark} If we think about the 2-crossed module $L\to M\to N$ as about a 2-groupoid with objects in $L$, 1-arrows in $L\times M$ and 2-arrows in $L\times N\times M$ then $L\to M\to N$-bundle gerbes should give an example of bigroupoid 2-torsors introduced in \cite{Bakovic}.

\subsection{Pullback} \label{Pull2}If $f :X\to X'$  then we put $f^*(\GG,\nn)=(f^*(\GG), f^*\nn)$; this will be again a 2-crossed module bundle gerbe.

\subsection{Definition}\label{Iso2CMBG}
We call two $(L\to M \to
N)$-bundle gerbes $(\GG,\nn)$ and $(\GG',\nn')$ over the same manifold $X$ stably
isomorphic if there exists a stable isomorphism $\qq\coloneqq (\QQ, \tilde \ll): \GG\to \GG'$ of  $(L\to M)$-bundle gerbes  such that $\nn'\mu(\qq) = \nn$  holds true for the induced isomorphism of trivial bundles $\mu(\qq): \mu(\GG)\to \mu(\GG')$. An $(L\to M\to N)$-bundle gerbe is trivial if it is stably isomorphic to the trivial  $(L\to M\to N)$-bundle gerbe $(((Y^{[2]}\times L, \partial_1{\rm pr}_L), Y, X, 1), {\rm pr}_N)$.

Pullbacks preserve stable isomorphisms, in particular a pullback of a trivial 2-crossed module bundle gerbe is again a trivial 2-crossed module bundle gerbe.

\subsection{Example}\label{NLT2}Let us notice that a general $(L\to M\to N)$-bundle gerbe is not necessarily locally trivial, although it is locally trivial as an $(L\to M)$-bundle gerbe. For a function $n: X\to N$ such that ${\rm Im}(n)$ is not a subset of ${\rm Im}(\partial_2)$ the $(L\to M\to N)$-bundle gerbe $(((Y^{[2]}\times L, \partial_1{\rm pr}_L), Y, X, 1), {\rm pr}_N.n{\rm pr}_X) $ is locally non-trivial. We will refer to such a 2-crossed module as 2-crossed module bundle gerbe defined by the $N$-valued function $n$ on $X$. Two such 2-crossed module bundle gerbes are stably isomorphic iff their respective functions $n$ and $n'$ are related by an $M$-valued function $m$ by $n'=\partial_2(m)n$. We will refer to such a stable isomorphism as being defined by the function $m$. Further, two such stable isomorphisms defined by respective functions $m$ and $m'$ are isomorphic iff they are related by an $L$-valued function $l$ on $X$ by $m'=\partial_1(l)m$.
\subsection{Example} \label{1toG3}Consider an $(1\to G\to N)$-bundle gerbe $(\GG, \nn)$. As a $(1\to G)$-bundle gerbe $\GG$ gives a principal $G$-bundle $P$ (more precisely a function $G$-valued function on $g$ on $Y^{[2]}$ satisfying on $Y^{[3]}$ the 1-cocycle relation). The trivializing section $\nn$ then gives an $N$ valued function $n$ on $Y$ such that $g_{12}n_2=n_1$ on $Y^{[2]}$ and hence, a trivialization of the left principal $G$-bundle $P$ under the map $G\to N$. Hence, a $(1\to G\to N)$-bundle gerbe is the same thing as a $(G\to N)$-bundle.

Obviously, isomorphic $(G\to N)$-bundles correspond to stably isomorphic $(1\to G\to N)$-bundle gerbes.

\subsection{Remark}\label{Indep2} Let $\wp: Y'\to X$ be an another surjective submersion and $f: Y'\to Y$ a map such that $\wp'=\wp f$. Then 2-crossed module bundle gerbes $(f^*\GG, \nn)$ and $(\GG,\nn)$ are stably isomorphic. This can be shown in a complete analogy to the case of a crossed module bundle gerbe (\ref{Indep1}).
\subsection{Change of the structure 2-crossed module}\label{Change3} If $(L\to M\to N) \to (L'\to M'\to N')$ is a morphism of 2-crossed modules, there is an obvious way to construct starting from an $(L\to M\to N)$-bundle gerbe $(\GG, \nn)$ an $(L'\to M'\to N')$-bundle gerbe $(\GG', \nn')$ by changing the structure crossed module of $\GG$ from $L\to M$ to $L'\to M'$ and putting $\nn' = \nu \nn$ where $\nu$ is the morphisms  $\nu:N\to N'$ entering the definition of the morphism of  two 2-crossed modules. Obviously, the change of the structure 2-crossed module preserves stable isomorphisms of 2-crossed module bundle gerbes.

\subsection{Definition}\label{IsoIso2CMBG} Let $(({\PP},Y,X,\ll),\nn)$ and $(({\PP'},Y',X,\ll'),\nn')$ be two 2-crossed module bundle gerbes and $(\QQ, \tilde\ll_\QQ)$ and $(\RR, \tilde\ll_\RR)$ two stable isomorphism between them; see (\ref{Iso2}). We call $(\QQ, \tilde\ll_\QQ)$ and $(\RR, \tilde\ll_\RR)$ isomorphic if there is an isomorphism $\underline\ll:\QQ \to \RR$ of crossed module bundles on $\bar Y=Y\times_X Y'$ such that (with an obvious abuse of notation) the diagram
\begin{equation}\label{iso2crossmodbunger}
\begin{CD}
 q^*\PP@>
 {\tilde\ll_\QQ}>> \QQ_1{q'}^*\PP'\QQ_2^{-1}\\
 @V {\rm id}
 VV @V  VV \hskip-0.5cm{\scriptstyle{\underline\ll_1}\underline\ll_2^{-1}}\\
q^*\PP@>
 {\tilde\ll_\RR}>>\RR_1{q'}^*\PP'\RR_2^{-1}
\end{CD}
\end{equation}
is commutative. Obviously, pullbacks preserve isomorphisms of stable isomorphisms.
\subsection{2-cocycles}\label{2-cocycle3}
Let $\wp:Y\to X$ be the surjective submersion, which was implicitly contained in the
above definition of a 2-crossed module bundle gerbe. Since also for 2-crossed module bundle gerbes it holds true that 2-crossed module bundle gerbes $(f^*\GG,\nn)$ and $(\GG,\nn)$ are stably isomorphic if the respective maps $\wp$ and $\wp'$ are related by a compatible map, we can again assume $Y=\coprod O_i$. For simplicity, we assume that the covering $\{O_i \}$ is a good one, in which case the $(L\to M \to N)$-bundle gerbe is
characterized
by transition functions $(n_i, m_{ij}, l_{ijk})$, $ n_i:O_{i}\to N$, $m_{ij}: O_{ij}\to M$ , $l_{ijk}:O_{ijk}\to L $ fulfilling 2-cocycle relations
$$n_i= \partial_2 (m_{ij})n_j$$
$$m_{ij} m_{jk} =\partial_1 (l_{ijk})m_{ik}$$
$$ l_{ijk}l_{ikl}=\,^{m_{ij}}l_{jkl}l_{ijl}$$ on $O_{ij}$, $O_{ijk}$ and $O_{ijkl}$, respectively.

In terms of 2-cocycles the stable isomorphism
$(l_{ijk}, m_{ij}, n_i) \sim (l_{ijk}', m_{ij}', n_i')$
is expressed by relations
$$n_i' = \partial_2 (m_i) n_i$$
$$m_{ij}'= m_i\partial_1 (l_{ij}) m_{ij} m_j^{-1}$$
$$\,^{m_i^{-1}}l_{ijk}' = l_{ij}\,^{m_{ij}} l_{jk}l_{ijk}
l_{ik}^{-1}$$

A trivial 2-crossed module bundle gerbe is described by transition functions
$$n_i = \partial_2 (m_i)$$
$$m_{ij}=m_i\partial_1(l_{ij}) m_j^{-1}$$ and $$^{m_i^{-1}}l_{ijk}=l_{ij}l_{jk}
 l_{ik}^{-1}$$

Locally, two collections of stable isomorphism data $(m_i, l_{ij})$ and $(m_i', l_{ij}')$ are isomorphic if
$$m_i' = \partial_1 (l_i)m_i $$
$$l_{ij}'= l_i l_{ij}l_{j}^{-1}$$

An $(L\to M\to N)$-bundle gerbe can be reconstructed from transition functions $(n_i, m_{ij}, l_{ijk})$ in a complete analogy with the case of an $(L\to M)$-bundle gerbe (\ref{2-cocycle}).
Again, it is easy to check that stable isomorphism classes of $(L\to M\to N)$-bundle gerbes are one to one with stable isomorphism classes of 2-cocycles. This follows rather simply from the corresponding statement for $(L\to M)$-bundle gerbes (\ref{2-cocycle}). Also, similarly to the case of crossed module bundles (\ref{2-cocycle}), we can consider a 2-category of $(L\to M\to N)$-bundle gerbes, with 1-arrows being stable isomorphisms and 2-arrows being isomorphism of stable automorphisms and similarly a 2-category of the 2-cocycles, but we will not use these.
\subsection{Lifting 2-crossed module bundle gerbe}\label{Lifting3}
Consider a Lie 2-crossed module $L\to M \to N$ such that $\ker(\partial_1)=1$ and
$\ker(\partial_2)={\rm Im}(\partial_1)$.
Put $G := L\backslash M $ and $Q := G\backslash N$. Assume that $L$ is a closed subgroup of $M$ and $G$ is a closed subgroup of $N$. So we have extensions of Lie groups
$$1\to L\stackrel{\partial_1}{\to} M \stackrel{\partial_2}{\to} N\stackrel{\pi_2}{\to} Q \to 1$$
$$1\to L\stackrel{\partial_1}{\to} M\stackrel{\pi_1}{\to}  G\to 1$$ and
$$1\to G\stackrel{\partial_2'}{\to} N\stackrel{\pi_2}{\to}  Q \to 1$$
such that $M\stackrel{\pi_1}{\to}G$ is an $(L\to M)$-bundle and $N\stackrel{\pi_1}{\to}Q$ is an $(G\to N)$-bundle.
Also, we have an exact sequence of pre-crossed modules
$$
\begin{CD}
 1 @>{}>>L@>
 {\partial_1}>> M @>\pi_1>> G @>{}>>1
 \\
 @V {} VV @V  VV @V  \hskip 2.5cm{\scriptstyle{\partial_2}}VV@V  \hskip 2.5cm{\scriptstyle{\partial_2'}}VV @V  VV\\
1 @>{}>>1 @> {}>> N @>{}>> N @>{}>>1
\end{CD}
$$
where $G$ is a normal subgroup of $N$ and also a morphisms of 2-crossed modules
$$
\begin{CD}
 L@>
 {\partial_1}>> M@>
 {\partial_2}>> N\\
 @V {}
 VV @V \pi_1 VV @V  \rm{id}VV\\
1@>
>>G@>
 {\partial_2'}>>N\,
\end{CD}
$$
Let us first notice that given a
$(G\to N)$-bundle $\PP=(P,\nn)$ on $X$, the left principal $G$-bundle $P$ can be lifted to an $(L\to M)$-bundle gerbe $\GG$ (\ref{Lifting2}), which will be actually an $(L\to M\to N$)-bundle gerbe. This is because of the identity $\partial_2'\pi_1 = \partial_2$ the trivialization $\nn$ of $\PP$ under $\partial_2'$ naturally defines a trivialization  of the principal $N$-bundle $\mu(\GG)$ (let us recall, $\mu: (L\stackrel{\partial_1}{\to}M)\to (1\to N)$ is a morphism of crossed modules). On the other hand, starting with an $(L\to M\to N)$-bundle gerbe $(\GG, \nn)$ with the 2-crossed module as above, we can change its structure 2-crossed module to $1\to G\to N$ in order to obtain a principal $(G\to N)$-bundle $\PP$.

Again, these two constructions are inverse to each other on sets of stable isomorphism classes of $(L\to M\to N)$-bundle gerbes and  isomorphism classes of $(G\to N)$-bundles (with $L$, $M$, $N$ and $G$ originating from a 2-crossed module as above). Given three $(G\to N)$-bundles $\PP, \PP'$ and $\PP''$ and isomorphisms $\PP\stackrel{f}{\to} \PP'$, $\PP'\stackrel{f'}{\to} \PP''$ and $\PP\stackrel{f''}{\to} \PP''$ such that $f'f=f''$ the corresponding lifting 2-crossed module bundle gerbes will be stably isomorphic, but the respective stable isomorphisms $\ff\ff'$ and $\ff''$ will be only isomorphic in general.

If we now consider a $Q$-valued function $q$ on $X$ we can lift it to a $(G\to N)$-bundle and after lifting this $(G\to N)$-bundle we obtain an $(L\to M\to N)$-bundle gerbe. Going into other direction, starting from an $(L\to M\to N)$-bundle gerbe and changing its structure 2-crossed module to $1\to G\to N$ we obtain a principal $G\to N$-bundle $\PP$, and further changing the $(G\to N)$ crossed module to $(1\to Q)$ we get a $Q$-valued function. The two above constructions are inverse to each other up to a stable isomorphism  of $(L\to M\to N)$-bundle gerbes.

The local description of the above constructions is also similar to the case of crossed module bundle gerbes (\ref{Lifting2}).
\subsection{Remark}\label{RightAction2} For an $(L\to M\to N)$-bundle gerbe $(\GG,\nn)=((P,\mm), Y, X, \ll),\nn)$ we recall from (\ref{1toG3}) that the trivializing section $\nn$ of the left principal $N$-bundle $\mu(\GG)$ defines an $N$-valued function $\nnn$ on $Y$ such that $\partial_2(\mm)=\nnn_1 \nnn_2^{-1}$. Let us recall that on the left principal $L$-bundle $P$ there is a compatible principal right $L$-action. Using the $N$-valued function $\nnn$ we can introduce a further principal right $L$-action on $P$, which will again commute with the principal left $L$-action. We will use the notation $(p,l)\mapsto p._\nnn l$ for $p\in P$, $l\in L$ for this principal right action of $L$ and put
$p._\nnn l\coloneqq p\,^{\nnn_2(y_1,y_2)}l$, where $p$ lies in the fibre over $(y_1,y_2)\in Y^{[2]}$ and $\nnn_2$ is the pullback to $Y^{[2]}$ of
$\nnn$ under the projection to the second factor of $Y^{[2]}$. It is easy to check that this formula indeed defines a principal right $L$-action
commuting with the principal left $L$-action on $P$.

Let now  $(\GG,\nn)=((P,\mm), Y, X, \ll),\nn)$ and $(\tilde\GG,\tilde\nn)=((\tilde P,\tilde
\mm), Y, X, \tilde \ll),\tilde \nn)$ be two 2-crossed module bundle gerbes. Let us again consider on $Y^{[2]}$ the Whitney sum $P\oplus \tilde P$
and introduce an equivalence relation on $P\oplus \tilde P$ by
$$(p._\nnn l, \tilde p)\sim_\nnn (p, l \tilde p)$$
and define $\bar P =P._\nnn \tilde P= P\oplus \tilde P/\sim_\nnn.$ We will denote an element of $P._\nnn \tilde P$ defined by equivalence class of $(p,\tilde p)\in P\oplus \tilde P$ as $[p,\tilde p]_\nnn$ in order to distinguish it from equivalence class $[p,\tilde p]\in P\tilde P$ defined previously in (\ref{mult1}).
Also, put
$$\bar\mm = \mm\,^{\nnn_2}\tilde\mm$$
It is easy to see that $\bar\PP\coloneqq(\bar P,\bar\mm)$ is an $(L\to M)$-bundle on $Y^{[2]}$. Let us notice that also $\partial_2(\bar\mm)=\bar \nnn_1\bar\nnn_2$ on $Y^{[2]}$ with
$$\bar\nn=\nn\tilde \nn$$
Now on $Y^{[3]}$ we do have the pullbacks $\PP_{12}$, $\tilde\PP_{12}$, $\bar\PP_{12}$, etc. An element of $\bar P_{12}\bar P_{23}$ is then given by
$((y_1,y_2,y_3), [[p, \tilde p]_\nnn , [p', \tilde p']_\nnn ])$ with $(y_1,y_2,y_3)\in Y^{[3]}$, $p\in P$ and $\tilde p\in \tilde P$ in the respective fibres of $P$ and $\tilde P$ over $(y_1,y_2)\in Y^{[2]}$, and $p'\in P$ and $\tilde p'\in \tilde P$  are in the respective fibres of $P$ and $\tilde P$ over $(y_2,y_3)\in Y^{[2]}$.
Finally, we define $\bar \ll:\bar P_{12}\bar P_{23} \to \bar P_{13}$ as
$$\bar \ll((y_1,y_2,y_3), [[p, \tilde p]_\nnn , [p', \tilde p']_\nnn ])\coloneqq ((y_1,y_2,y_3), [\ll([p, p'], \tilde\ll[\tilde p, \tilde p']]_\nnn )$$
Now it is a rather lengthy but a straightforward check to establish the following proposition.
\subsection{Proposition} \label{mult2}$(\bar \GG,\bar\nn)\coloneqq((\bar P, \bar \mmm ), Y, X, \bar \ll),\bar \nn)$ defines an $(L\to M\to N)$-bundle gerbe, the product of $(L\to M\to N)$-bundle gerbes $(\GG,\nn)=((P,\mm), Y, X, \ll),\nn)$ and $(\tilde\GG,\tilde\nn)=((\tilde P,\tilde \mm), Y, X, \tilde \ll),\tilde \nn)$.
\subsection{Example}\label{ProdLNT2} If $(\GG,\nn)=(((P=Y^{[2]}\times L, \partial_1{\rm pr}_L), Y, X, 1), {\rm pr}_N.n{\rm pr}_X)$ and $(\tilde \GG,\tilde \nn)=(((\tilde PY^{[2]}\times L, \partial_1{\rm pr}_L), Y, X, ), {\rm pr}_N.\tilde n{\rm pr}_X)$ are two $(L\to M\to N)$-bundle gerbes defined by two respective $N$-valued functions $n$ and $\tilde n$ on $X$ (\ref{NLT2}) then their product is explicitly described again as an $(L\to M\to N)$-bundle gerbe defined by the function $n\tilde n$ by identifying $[(y_1,y_2,l), (y_1,y_2,\tilde l)]\in PP'$ with $(y_1,y_2, l \,^{n(x)}\tilde l)\in Y^{[2]}\times L)$. Here $(y_1,y_2)\in Y^{[2]}$ live in the fibre over $x\in X$.
\subsection{Remark}\label{Group3}
The above product defines a groups structure on stable isomorphism classes of $(L\to M \to N)$-bundle gerbes. The unit is given by the class of the trivial $(L\to M \to N)$-bundle gerbe  $(((Y^{[2]}\times L, \partial_1{\rm pr}_L), Y, X, 1), {\rm pr}_N)$. We will give an explicit (local) formula for the inverse later.
Let us notice that the relation between the stable isomorphism classes of lifting $(L\to M\to N)$-bundle gerbes described above (\ref{Lifting3}) and $Q$-valued functions (and stable isomorphism classes of $(G \to N)$-bundle gerbes) is compatible with the respective multiplications.

\subsection{Product on 2-cocycles}\label{LocalProd}
The product formulas
for the corresponding transition functions (2-cocycles) of the product
$\bar \GG=\GG \tilde\GG$ of two 2-crossed module bundles are given by
$$\bar n_i = n_i \tilde n_i$$
$$\bar m_{ij}=m_{ij}\,^{n_j}\tilde m_{ij}$$
$$\bar l_{ijk} =l_{ijk}\,^{m_{ik}}\{{m_{jk}}^{-1},
\,^{n_j}\tilde m_{ij}\}\,^{n_i}\tilde l_{ijk}$$
The inverse $(n_i, m_{ij}, l_{ijk})^{-1}$ is given by
$$(n_i^{-1}, \,^{n_j^{-1}}m_{ij}^{-1},\,^{n_k^{-1}}\{m_{jk}^{-1},m_{ij}^{-1}\}^{-1} \,^{n_i^{-1}}l^{-1}_{ijk})$$
\subsection{Remark} \label{lax2}Let us forget, for the moment, about the "horizontal" composition in the Gray  3-groupoid corresponding to the 2-crossed module $L\to M\to N$. The two "vertical" compositions define a strict 2-groupoid (a strict topological 2-category), which we will denote as $2\CC$. Let us again consider the topological $\OO$ category defined by the good covering $\{O_i\}$. A 2-cocycle is the same thing as a continuous, normal pseudo-functor from $\OO$ to $2\CC$.
Now we can use the fact that the horizontal composition in a topological Gray 3-category defines a continuous cubical functor $\mathfrak F:2\CC \times 2\CC \to 2\CC$ from the cartesian product $2\CC \times 2\CC$ to $2\CC$ \cite{DayStreet}. We may use the following property of cubical functors, which follows almost immediately from definition. If  $\FF$ and $\GG$  are two continuous normal pseudo-functors from $\OO$ to $2\CC$ then $\mathfrak F(\FF,\GG)$ is a pseudo-functor from $\OO$ to $2\CC$. Hence, we obtain a product on 2-cocycles, which is the same as the one given above (\ref{mult2}).
\section{2-crossed module bundle 2-gerbes}\label{2CM2BG}

Consider again a surjective submersion $\wp:Y\to X$. Let, as before, $p_{ij}:
Y^{[3]} \to Y^{[2]}$ denotes the projection to the $i$-th and
$j$-th component, and similarly for projections of higher fibred powers $Y^{[n]}$ of
$Y$. Let $L\stackrel{\partial_1}{\to} M \stackrel{\partial_2}{\to} N$ be a 2-crossed module.
\subsection{Definition} A 2-crossed module bundle 2-gerbe is defined by a quintuple $(\mathfrak G, Y, X, \mm, \ll)$, where $\mathfrak G=(\GG,\nn)$
is a 2-crossed module bundle gerbe  over $Y^{[2]}$,
$$\mm: \mathfrak G_{12}\mathfrak G_{23}\to \mathfrak G_{13}$$
is a stable isomorphism  on $Y^{[3]}$ of the product $\mathfrak G_{12}\mathfrak G_{23}$
of the pullback 2-crossed module bundle gerbes $\mathfrak G_{12}=p_{12}^*\mathfrak G$ and $\mathfrak G_{23}=p_{23}^*\mathfrak G$ and the pullback 2-crossed module bundle gerbe $\mathfrak G_{13}=p_{13}^*\mathfrak G$, and
$$\ll: \mm_{124}\mm_{234}\to \mm_{134}\mm_{123}$$
is an isomorphism of the composition of pullbacks of stable isomorphisms $p_{124}^*\mm$ and $p_{234}^*\mm$ and the composition of pullbacks of stable isomorphisms $p_{123}^*\mm$ and $p_{134}^*\mm$ on $Y^{[4]}$.
On $Y^{[5]}$, the isomorphism $\ll$ should satisfy the obvious
coherence relation
$$\ll_{1345}\ll_{1235} = \ll_{1234}\ll_{1245}\ll_{2345}.$$
\subsection{Abelian bundle 2-gerbes}\label{Abelian2}Abelian bundle 2-gerbes as introduced in \cite{Stev2-gerbes}, \cite{StevPhD}, \cite{CareyMuWHighrBG} are $(U(1)\to 1\to 1)$-bundle 2-gerbes. If $A\to 1\to 1$ is a 2-crossed module then $A$ is necessarily
an abelian group and an abelian bundle 2-gerbe can be identified as an $(A\to 1 \to 1)$-bundle 2-gerbe.
\subsection{Example} \label{Example4}Consider an $(1\to G\to N)$-bundle 2-gerbe $(\mathfrak G,Y,X,\mm,\ll)$. The $(1\to G\to N)$-bundle gerbe on $Y^{[2]}$ gives a $(G\to N)$-bundle $\PP$ on $Y^{[2]}$. The stable isomorphism $\mm:\mathfrak G_{12} \mathfrak G_{23} \to  \mathfrak G_{13}$ gives on $Y^{[3]}$ an isomorphism $\gg:\PP_{12} \PP_{23} \to  \PP_{13}$ satisfying on $Y^{[4]}$ the associativity condition $\gg_{124}\gg_{234}=\gg_{134}\gg_{123}$ since the first Lie group of the 2-crossed module  $(1\to G\to N)$ is trivial.  Hence, a $(1\to G\to N)$-bundle 2-gerbe is the same thing as a $(G\to N)$-bundle gerbe. Obviously, stably isomorphic $(1\to G\to N)$-bundle 2-gerbes correspond to stably isomorphic $(G\to N)$-bundle gerbes.
\subsection{Pullback}\label{Pull4} If $f :X\to X'$ is a map then we can pullback $Y\to X$ to $f^*(Y)\to X'$ with a map $\tilde f: f^*(Y)\to Y$ covering $f$.
There are induced maps $\tilde f^{[n]} : f^*(Y )^{[n]} \to Y^{[n]}$.
The pullback $f^*({\mathfrak G},Y,X,\mm, \ll)\coloneqq (\tilde f^{[2]*}{\mathfrak G },f^*(Y),f(X),\tilde f^{[3]*}\mm, \tilde f^{[4]*}\ll)$ is again an $(L\to M\to N)$-bundle 2-gerbe.
\subsection{Definition}\label{Iso$}
Two 2-crossed module bundle 2-gerbes $((\mathfrak G,Y,X,\mm,\ll)$ and $(\mathfrak G',Y',X,\mm', \ll')$ are
stably isomorphic if there exists a 2-crossed module bundle gerbe $\mathfrak Q \to \bar Y=Y\times_X Y'$ such that over $\bar Y^{[2]}$ the 2-crossed module bundle gerbes $q^*\mathfrak G$ and
$\mathfrak Q_1 q'^*\mathfrak G'\mathfrak Q_2^{-1}$ are stably isomorphic. Let $\tilde\mm$ be the stable isomorphism
$\tilde\mm: q^*\mathfrak G  \to \mathfrak Q_1 q'^*\mathfrak G'\mathfrak Q_2^{-1}$. Then we ask  on $Y^{[3]}$ (with an obvious abuse of notation) an existence of an isomorphism $\tilde \ll$ of stable isomorphisms
$$\tilde \ll:\mm'\tilde\mm_{23}\tilde\mm_{12}\to \tilde\mm_{13}\mm $$ fulfilling on $Y^{[4]}$
$$\ll_{1234}\tilde\ll_{124}\tilde\ll_{234}= \tilde\ll_{134}\tilde\ll_{123}\ll'_{1234}$$
Here $q$ and $q'$ are projections onto first and second factor of $\bar Y=Y\times_X
Y'$ and $\mathfrak Q_1$ and $\mathfrak Q_2$ are the pullbacks of  $\mathfrak Q$ to $\bar Y^{[2]}$ under respective
projections $p_1$, $p_2$ form $\bar Y^{[2]}$ to $\bar Y$, etc.

A 2-crossed module bundle 2-gerbe $(\mathfrak G,Y,X,\mm,\ll)$ is called trivial if it is stably isomorphic to the trivial 2-crossed module bundle 2-gerbe $({\mathcal T}, Y, X, 1, 1)$, where ${\mathcal T}$ is the trivial 2-crossed module bundle gerbe $(((Z^{[2]}\times L, \partial_1{\rm pr}_L), Z, Y^{[2]}, 1), {\rm pr}_N))$. Pullbacks preserve stable isomorphisms, a pullback of a trivial 2- crossed module bundle 2-gerbe is again a trivial 2-crossed module bundle 2-gerbe.

If $Y=X$ then the 2-crossed module bundle 2-gerbe is trivial.
\subsection{Definition}\label{IsoIso4} Let $({\mathfrak G},Y,X,\mm, \ll)$ and $({\mathfrak G'},Y',X,\mm',\ll')$ be two 2-crossed module bundle 2-gerbes and $(\mathfrak Q, \tilde\mm_{\mathfrak Q}, \tilde\ll_{\mathfrak Q} )$ and $(\mathfrak R, \tilde\mm_{\mathfrak R}, \tilde\ll_{\mathfrak R} )$ two stable isomorphism between them. We call these two stable isomorphisms stably isomorphic if there is a stable isomorphism $\underline\mm:\mathfrak Q \to \mathfrak R$ of 2-crossed module bundles on $\bar Y=Y\times_X Y'$ such that (with an obvious abuse of notation) the diagram
$$
\begin{CD}
 q^*{\mathfrak G}@>
 {\tilde\mm_{\mathfrak Q}}>> {\mathfrak Q}_1{q'}^*{\mathfrak G}'{\mathfrak Q}_2^{-1}\\
 @V \rm{id}
 VV @V  VV \hskip-0.5cm{\scriptstyle{\underline\mm_1}{\underline\mm_2}^{-1}}\\
q^*{\mathfrak G}@>
 {\tilde\mm_{\mathfrak R}}>>{\mathfrak R}_1{q'}^*{\mathfrak G'}{\mathfrak R}_2^{-1}
\end{CD}
$$
commutes up to an isomorphism of stable isomorphisms
$$\underline \ll:\tilde\mm_{\mathfrak Q}{\underline\mm_1}{\underline\mm_2}^{-1}\to \tilde\mm_{\mathfrak R}$$
on $\bar Y^{[2]}$, fulfilling on $\bar Y^{[3]}$
$$\tilde \ll_{\mathfrak Q}\underline \ll_{13}=\underline \ll_{12}\underline \ll_{23}\tilde\ll_{\mathfrak R}$$
\subsection{Remark} \label{Indep4}
Let $\wp: Y'\to X$ be an another surjective submersion and $f: Y'\to Y$ a map such that $\wp'=\wp f$. Then 2-crossed module bundle 2-gerbes $(f^*\mathfrak G, Y', X, f^{[3]*}\mm, f^{[4]*}\ll)$ and $(\mathfrak G, Y, X, \mm,\ll)$ are stably isomorphic. It follows that locally each 2-crossed module bundle gerbe $\GG$ is trivial. The arguments to show the above two statements are completely analogous to the case of a crossed module bundle gerbe (\ref{Indep2}).
\subsection{Change of the structure 2-crossed module}\label{Change4} If $(L\to M\to N) \to (L'\to M'\to N')$ is a morphism of crossed modules, there is an obvious way to construct, starting from an $(L\to M\to N)$-bundle 2-gerbe $(\mathfrak G,Y,X,\mm,\ll)$, an $(L'\to M'\to N')$-bundle 2-gerbe $(\mathfrak G',Y,X,\mm',\ll')$  by changing the structure 2-crossed module of $\mathfrak G$ from $L\to M \to N$ to $L'\to M'\to N'$.

\subsection{3-cocycles}\label{3-cocycle}
Let $\wp:Y\to X$ be the surjective submersion, which was implicitly contained in the above definition of a 2-crossed module bundle 2-gerbe. Let us
recall (\ref{Indep4}) that also for 2-crossed module bundle 2-gerbes it holds true that 2-crossed module bundle 2-gerbes $(f^*\mathfrak G, Y', X, f^{[3]*}\mm, f^{[4]*}\ll)$ and $(\mathfrak G, Y, X, \mm,\ll)$ are stably isomorphic if the respective maps $\wp$ and $\wp'$ are related by a compatible map $f$. Hence, we can again assume $Y=\coprod O_i$. For simplicity, we assume that the covering $\{O_i \}$ is a good one, in which case the $(L\to M \to N)$-bundle gerbe can be described by transition functions
$(n_{ij}, m_{ijk},l_{ijkl} )$ $n_{ij}: O_{ij}\to N$, $m_{ijk}:
O_{ijk}\to M$ and $l_{ijkl}: O_{ijkl}\to L$ satisfying
\begin{equation}\label{3cocycle}
\begin{array}{c}
n_{ij}n_{jk}= \partial_2(m_{ijk}) n_{ik}\\
m_{ijk}m_{ikl}= \partial_1(l_{ijkl})\,^{n_{ij}}
m_{jkl}m_{ijl}\\
l_{ijkl}\,^{\,^{n_{ij}}m_{jkl}}(l_{ijlm}) \,^{n_{ij}}l_{jklm}= \,^{m_{ijk}}l_{iklm}\{m_{ijk},\,^{n_{ik}} m_{klm}\}\,^{\,^{n_{ij}n_{jk}}m_{klm}}
(l_{ijkm})
\end{array}
\end{equation}
We shall not give explicit formulas for stably isomorphic transition functions and for transition functions of a trivial 2-crossed module bundle 2-gerbe, since
the respective expressions are rather cumbersome and not very illuminating.

Now we briefly describe how an $(L\to M)$-bundle gerbe can be
reconstructed from transition functions $(n_{ij}, m_{ijk}, l_{ijkl})$. Put $Y=\coprod O_i$. On each nonempty $O_{ij}$ consider the $(L\to M\to
N)$-bundle gerbe ${\mathfrak{G}}_{ij}$ defined by the function $n_{ij}:O_{ij}\to N$ as in (\ref{NLT2}). Hence, on $Y^{[2]}$ we have the $(L\to
M\to N)$-bundle gerbe given by ${\mathfrak{G}}=\coprod_{ij}{\mathfrak{G}}_{ij}$. Now, we recall the explicit descriptions of the multiplication
(\ref{ProdLNT2}) and stable isomorphisms (\ref{NLT2}) of two $(L\to M\to M)$-bundle gerbes defined by their respective $N$-valued functions. Also,
recall the description of isomorphisms of stable isomorphism in case of such $(L\to M\to M)$-bundle gerbes. Using the 3-cocycle relations, it is
now straightforward to show that the collection of functions $m_{ijk}$ defines a stable isomorphism of ${\mathfrak{G}}_{12}{\mathfrak{G}}_{23}$
and ${\mathfrak{G}}_{13}$ on $Y^{[3]}$ satisfying on $Y^{[4]}$ the associativity condition up to the an isomorphism defined by the collection of
functions $l_{ijkl}$, which fulfils the coherence relation on $Y^{[5]}$.

It might be interesting to examine possible 3-categorial aspects of the above constructions.
\subsection{Lifting 2-crossed module bundle 2-gerbe}\label{Lifting4}
Consider again a Lie 2-crossed module $L\to M \to N$ such that $\ker(\partial_1)=1$ and
$\ker(\partial_2)=\rm{Im}(\partial_1)$.
Put $ G :=L\backslash M$, $Q:= G\backslash N$ and consider the extensions of Lie groups
$$1\to L\stackrel{\partial_1}{\to} M \stackrel{\partial_2}{\to} N\stackrel{\pi_2}{\to} Q \to 1$$
$$1\to L\stackrel{\partial_1}{\to} M\stackrel{\pi_1}{\to}  G\to 1$$
$$1\to G\stackrel{\partial_2'}{\to} N\stackrel{\pi_2}{\to}  Q \to 1$$
We have the exact sequence of pre-crossed modules
$$
\begin{CD}
 1 @>{}>>L@>
 {\partial_1}>> M @>\pi_1>> G @>{}>>1
 \\
 @V {} VV @V  VV @V  \hskip 2.5cm{\scriptstyle{\partial_2}}VV@V  \hskip 2.5cm{\scriptstyle{\partial_2'}}VV @V  VV\\
1 @>{}>>1 @> {}>> N @>{}>> N @>{}>>1
\end{CD}
$$
where $G$ is a normal subgroup of $N$ and the
 morphisms of 2-crossed modules
$$
\begin{CD}
 L@>
 {\partial_1}>> M@>
 {\partial_2}>> N\\
 @V {}
 VV @V \pi_1 VV @V  \rm{id} VV\\
1@>
>>G@>
 {\partial_2'}>>N\,
\end{CD}
$$
Assume again that $L$ is a closed subgroup of $M$ and $G$ is a closed subgroup of $N$. Let $P$ be a left principal $Q$-bundle over $X$. Let us consider the corresponding lifting $(G\to N)$-bundle gerbe. This in particular means that on $P^{[2]}$ we have a $(G\to N)$-bundle which can be lifted to an $(L\to M\to N)$-bundle gerbe $\mathfrak G$ on $P^{[2]}$ (\ref{Lifting3}). It follows that the 2-crossed module bundle gerbes $\mathfrak G_{12}\mathfrak G_{23}$ and $\mathfrak G_{13}$ are stably isomorphic with a stable isomorphism $\mm$. This follows from the above mentioned fact that in case of Lie groups $L$, $M$, $N$ and $Q$ as above stable isomorphisms classes of $(L\to M\to N)$-bundle gerbes are one to one to $Q$-valued functions and that this correspondence respects the respective multiplications. Obviously, such a stable isomorphism $\mm$  fulfills the associativity condition on $Y^{[4]}$ up to an isomorphism fulfilling the coherence relation on $Y^{[5]}$ because of $\ker\partial_1=1$.

Going in the opposite direction, let us consider an $(L\to M \to N)$-bundle 2-gerbe $(\mathfrak G,X,Y,\mm, \ll)$ with the 2-crossed module $(L\to M \to N)$ as above. Changing the structure 2-crossed module to $1\to G\to N$, we obtain a $G\to N$-bundle gerbe $(\GG,
n)$ on $X$. After  changing its structure crossed module to $1\to Q$ we obtain a left principal $Q$-bundle on $X$.

The two above constructions are inverse to each other on the sets od stable isomorphisms of $(L\to M\to N)$-bundle 2-gerbes and isomorphism classes of principal $Q$-bundles (with $(L\to M\to N)$ and $Q$ as above).
\subsection{Remark}\label{lifting6}Obviously, one can reinterpret the above lifting construction (\ref{Lifting4}) as follows. From the 3-term exact sequence
$1\to L\stackrel{\partial_1}{\to} M\stackrel{\pi_1}{\to}  G\to 1$
and the right principal $(G\to N)$-bundle $N\to Q$ (given by the 3-term exact sequence $1\to G\stackrel{\partial_2'}{\to} N\stackrel{\pi_2}{\to}  Q \to 1$) we can construct a lifting $(L\to M)$-bundle gerbe on $Q$. This lifting bundle gerbe will actually be an $(L\to M\to N)$-bundle gerbe $\mathfrak G$ (\ref{Lifting3}). If now, as above, $P$ is a left principal $Q$-bundle over $X$ then we can use the corresponding division map $P^{[2]}\to G$ to pullback the 2-crossed module gerbe $\mathfrak G$ from $G$ to $P^{[2]}$. What we get is obviously a 2-crossed module bundle 2-gerbe stably isomorphic to the lifting bundle 2-gerbe of (\ref{Lifting4}).
\subsection{Twisting crossed module bundle gerbes with abelian bundle 2-gerbes} \label{Abelian4} Twisted crossed module bundle gerbes as discussed here were introduced in \cite{AschJur}. A more general concept of twisting has been introduced recently in \cite{SSS1}.

Let us consider a crossed module  $L\stackrel{\delta}{\to}M$. It follows, $\ker \delta\eqqcolon A$ is abelian. Putting $Q\coloneqq {\rm{coker}}\delta$ we have an exact sequence
$$0\to A\stackrel{\partial}{\to} L\stackrel{\delta}{\to}M \to Q\to 1$$
Hence, as in the above example, we have a 2-crossed module $$A\stackrel{\partial}{\to} L\stackrel{\delta}{\to}M$$
with $\ker\partial=0$ and $\rm{Im}\partial=\ker\delta$.
Consider, as above, extensions of Lie groups
$$1\to A\stackrel{\partial}{\to} L\stackrel{\pi_1}{\to}  G\to 1$$ and
$$1\to G\stackrel{\delta'}{\to} M\stackrel{\pi_2}{\to}  Q \to 1$$
However, now we have an exact sequence of crossed modules
$$
\begin{CD}
 1 @>{}>>A@>
 {\partial}>> L @>\pi_1>> G @>{}>>1
 \\
 @V {} VV @V  VV @V  \hskip 2.5cm{\scriptstyle{\delta}}VV@V  \hskip 2.5cm{\scriptstyle{\delta'}}VV @V  VV\\
1 @>{}>>1 @> {}>> M @>{}>> M@>{}>>1
\end{CD}
$$
where $G$ is a normal subgroup of $M$. As before, we do have a
 morphisms of 2-crossed modules
$$
\begin{CD}
 A@>
 {\partial}>> L@>
 {\delta}>> M\\
 @V {}
 VV @V \pi_1 VV @V \rm{id} VV\\
1@>
>>G@>
 {\delta'}>>M\,
\end{CD}
$$
Assume again that $A$ is a closed subgroup of $L$ and $G$ is a closed subgroup of $M$. This means that staring from a crossed module
$L\stackrel{\delta}{\to}M$ and a principal $Q$-bundle $P$ we can construct a lifting 2-crossed module bundle 2-gerbe with abelian $A=\ker \partial$.
Let us further assume that what we have here is a central extension of $L$ by $A$, and that $M$ acts trivially on $A$. \label{central}
We may choose, in order to be more concrete, $A=U(1)$. Let us assume that the lifting bundle 2-gerbe $\mathrm{G})$ is
described locally, with respect to a good covering, by a 3-cocycle $(m_{ij}, l_{ijk}, a_{ijkl})$ \begin{equation}\label{bun2gerG}
\begin{array}{c}
m_{ij}m_{jk}= \delta(l_{ijk}) m_{ik}\\
l_{ijk}l_{ikl}= \partial(a_{ijkl})\,^{m_{ij}} l_{jkl}l_{ijl}\\
a_{ijkl}a_{ijlm} a_{jklm}= a_{iklm}a_{ijkm}.
\end{array}
\end{equation}
The collection of $U(1)$-valued functions $a_{ijkl}$ on the quadruple intersections represents a \v Cech class in $H^3(X, U(1))$ or
correspondingly a class in $H^4(X,\BZ)$. We may think of it as representing an abelian bundle 2-gerbe $\mathrm{A}$. If we assume that
$\mathrm{A}$ is trivial, we have
$$a_{ijkl}=\tilde a_{ijk}\tilde a_{ikl}\tilde a_{jkl}^{-1}\tilde a_{ijl}^{-1}.$$
Also, we see that we have a 2-cocycle $(l_{ijk}\partial(a_{ijk})^{-1}, m_{ij})$ representing a possibly non-trivial $(L\to M)$-bundle gerbe $\GG$. Obviously, the $(U(1)\to 1 \to 1)$-bundle 2-gerbes represented by non-trivial classes in $H^3(X, U(1))$ represent obstructions to lift a $(G\to M)$-bundle gerbe (and hence also a $Q$-bundle) to an $(L\to M)$-bundle gerbe.
Further, if $\tilde a_{ijk}$ and $\tilde a'_{ijk}$ represent two trivialization of $a_{ijkl}$ then $\tilde a_{ijk}(\tilde a'_{ijk})^{-1}$ represents a \v Cech class in $H^2(X, U(1))$ or correspondingly a class in $H^3(X,\BZ)$. We may think of it as representing an abelian bundle gerbe, i.e, the $(U(1)\to 1)$-bundle gerbe, $\mathcal{A}$.
We can summarize the above discussion in the following proposition.
\subsection{Proposition} In the notation of (\ref{Abelian4}):

i) A  principal $Q$-bundle on $X$ can be lifted to an $(L\to M)$-bundle gerbe if and only if the corresponding obstruction $(A\to 1\to 1)$-bundle 2-gerbe $\mathrm{A}$ is trivial.

ii) If non-empty, the set of stable isomorphism classes of those $(L\to M)$-bundle gerbes, which are liftings of $Q$-principal bundles  from the same isomorphism class, is freely and transitively acted on by the group of stable isomorphism classes of $(A\to 1)$-bundle gerbes.

These statements remain true in case when the principal $Q$-bundles and their isomorphisms classes are replaced by  $(G\to M)$-bundle gerbes and their stable isomorphism classes.
\subsection{Remark} Of course, the above lifting always exists when the 4-term exact sequence $1\to A\to L\to M \to Q \to 1$ corresponds to a trivial class in $H^3(Q,A)$ \cite{MacLane},\cite{Brown}, the third $Q$-cohomology with values in $A$.  The above lifting also trivially exists when  $X$ doesn't admit nontrivial $(A\to 1\to 1)$-bundle 2-gerbe, i.e., when $[X, B^2 A]$ is trivial.
\subsection{A remark on string structures} \label{StringStr} Let $Q$ be a simply-connected compact simple Lie group. Associated to $Q$ there is a crossed module $L\to M$ of infinite dimensional Fr\'echet Lie groups with $L:=\widehat {\Omega Q}$ and $ M:=P_0 Q$, where $\widehat {\Omega Q}$ is centrally extended group of based smooth loops in $Q$ and $P_0 Q$ is the group of smooth paths in $Q$ that start at the identity \cite{BaezCransSchrStev}. Hence in the notation of (\ref{Abelian4}) we do have $A=U(1)$, and $G=\Omega Q$. Let us notice (see \cite{StolzTeich}, \cite{BaezCransSchrStev}, \cite{Henriq}) that, in the situation as above (\ref{central}), the classifying space $BU(1)=K(\BZ, 2)$ can be equipped with a proper group structure and a topological group $String(Q)$ can be defined fitting an exact sequence of groups
$1\to K(\BZ,2)\to String(Q)\to Q \to 1$. Also, it is known \cite{Jurco}, \cite{BaezStevClass} that the categories of $(L\to M)$-bundle gerbes  and principal $String(Q)$-bundles are equivalent. A string structure is, by definition, a lift of a principal $Q$-bundle to a principal $String(Q)$-bundle and hence equivalently a lift of a $(G\to M)$-bundle gerbe to an $(L\to M)$-bundle gerbe. Thus, the above discussion applies to the existence of string structures and their classification as well.
\subsection{Remark}
A crossed square $(L\to A)\to (B\to N)$  \cite{LodayCrSq} of Lie groups gives a 2-crossed module, namely $L\to A\rtimes B \to N$ (see, e.g., \cite{MutPor}). A definition of a crossed square bundle 2-gerbe could possibly be read of from \cite{BreenSchr}, \cite{BreenNotes} \cite{BreenAst}. It would be interesting to compare these bundle 2-gerbes with $L\to A\rtimes B \to N$-bundle 2-gerbes defined in this paper.

\end{document}